\documentclass[12pt]{amsart}
\usepackage{amsmath}
\usepackage{latexsym}

\theoremstyle{remark}

\theoremstyle{plain}

\newtheorem{lemma}{Lemma}[section]
\newtheorem{theorem}{Theorem}

\numberwithin{equation}{section}
\DeclareMathAlphabet{\curly}{U}{rsfs}{m}{n}

\theoremstyle{definition}

\newtheorem*{Acknowledgement}{Acknowledgement} 

\newcommand{\lam}{\ensuremath{\lambda}}

\renewcommand{\a}{\ensuremath{\alpha}}
\renewcommand{\b}{\ensuremath{\beta}}
\newcommand{\del}{\ensuremath{\delta}}
\newcommand{\eps}{\ensuremath{\varepsilon}}
\newcommand{\sql}[1]{\ensuremath{\sqrt{\log{#1}}}}
\def\RR{\mathbb{R}}

\newcommand{\e}{\varepsilon}
\newcommand{\FF}{\curly{F}}

\newcommand{\vv}{\ensuremath{\mathbf{v}}}
\newcommand{\bE}{\ensuremath{\mathbf{E}}}
\newcommand{\uu}{\ensuremath{\mathbf{u}}}
\newcommand{\ww}{\ensuremath{\mathbf{w}}}

\newcommand{\zz}{\ensuremath{\mathbf{z}}}
\newcommand{\ba}{{\ensuremath{\boldsymbol{\alpha}}}}
\newcommand{\bb}{{\ensuremath{\boldsymbol{\beta}}}}
\newcommand{\bmu}{{\ensuremath{\boldsymbol{\mu}}}}
\newcommand{\brho}{{\ensuremath{\boldsymbol{\rho}}}}
\newcommand{\pfrac}[2]{\left(\frac{#1}{#2}\right)}
\newcommand{\ds}{\displaystyle}
\def\({\left(}
\def\){\right)}
\newcommand{\perm}{\operatorname{perm}}
\newcommand{\one}{\ensuremath{\mathbf{1}}}

\newcommand{\be}{\begin{equation}}
\newcommand{\ee}{\end{equation}}
\newcommand{\benn}{\begin{equation*}}   
\newcommand{\eenn}{\end{equation*}}

\begin{document} 

\title{On Bombieri's asymptotic sieve} 
\author{Kevin Ford}

\thanks{Research supported by National Science Foundation grants
DMS-0070618 and DMS-0301083.}

\address{Department of Mathematics, University of Illinois at 
Urbana-Champaign, Urbana, IL 61801  USA} 
 
\begin{abstract} 
If a sequence $(a_n)$ of non-negative real numbers has
``best possible'' distribution in arithmetic progressions,
Bombieri showed that one can deduce an asymptotic formula for the sum
$\sum_{n\le x} a_n \Lambda_k(n)$ for $k\ge 2$.  By constructing 
appropriate sequences, we show that any weakening of the well-distribution
 property is not sufficient to deduce the same conclusion.
\end{abstract}

\maketitle 

\section{Introduction}\label{sec:intro}               %

Many of the most famous problems in number theory
can be described in terms of estimating the number of
primes in an integer sequence.  More generally, given a sequence
$(a_n)$ of positive real numbers (e.g. the characteristic function
of a set of natural numbers), one can ask for bounds on the sum
$$
S_1(x) = \sum_{n\le x} a_n \Lambda(n),
$$
where $\Lambda$ is the von Mangoldt function. 
Removing from the sequence those terms with $n$ divisible by
a prime $\le z$ leaves behind only terms with $n$ composed of at most
$\lfloor \frac{\log x}{\log z}\rfloor$ prime factors.  If $z>\sqrt{x}$ then
only terms with $n$ prime are left.  Motivated by this simple fact,
the modern {\it sieve} was created by V. Brun (\cite{Br1}, \cite{Br2})
to attack such
problems, in particular the Twin Prime Conjecture and Goldbach's 
Conjecture.  Estimating the number of ``unsifted'' elements is usually
accomplished by means of a weighted form of inclusion-exclusion,
its precision entirely determined by the regularity of the sequence on the
arithmetic progressions $0\mod d$ for squarefree $d$ (see the monographs
\cite{Gre} and \cite{HR} for more about sieve procedures).  Writing
$$
A_d(x) = \sum_{\substack{n\le x \\ d|n }} a_n,
$$
one postulates the existence of a multiplicative function
$g$ so that 
$$
A_d(x) = g(d) A(x) + r_d(x),
$$
where $A(x)$ is an approximation to $A_1(x)$ and the ``remainders'' $r_d(x)$
are small in some average sense.   A typical hypothesis is
$$
R(\nu): \quad \forall B>0, \; 
\sum_{d\le x^{\nu}} |r_d(x)| \ll_{\nu,B} 
\frac{A(x)}{\log^B x}. 
$$

One also needs mild growth conditions on $A(x)$ and regularity conditions on
$g$.  There is some flexibility in choosing these conditions (see e.g.
\cite{Bo2}, \cite{FI1}, \cite{FI2}, \cite{Gre}, \cite{HR}), and generally 
these are easy to verify in practice.
We say that a sieve problem has {\it sifting density} or {\it dimension} 
$\kappa$ if $g(p)$ is about $\kappa/p$ on average over primes $p$.
In the important special case $\kappa=1$, one expects for many problems that 
\be\label{S1asym}
S_1(x) \sim H A(x), \qquad H = \prod_p (1-g(p)) (1-1/p)^{-1}.
\ee
For example, for the twin 
prime problem, we take $a_n=\Lambda(n+2)$, $A(x)=x$,
$g(d)=\frac{1}{\phi(d)}$ for
odd $d$ and $g(d)=0$ for even $d$.  It is known that $R(\nu)$ holds
for all $\nu < 1/2$ (the Bombieri-Vinogradov theorem), and it
is conjectured that $R(\nu)$ holds for all $\nu<1$.

That sieve methods cannot produce \eqref{S1asym} was discovered by
Selberg \cite{Se} in the 1940s.  His example is $a_n=1+\lambda(n)$, where
$\lambda(n)=(-1)^k$ if $n$ is the product of $k$ primes (not necessarily
distinct).  With $A(x)=x$ and $g(d)=1/d$, $R(\nu)$ holds for all $\nu<1$,
but $a_n=0$ for prime $n$ and
$$
S_1(x) = O(\sqrt{x}) = O(A(x) x^{-1/2}).
$$
In a sense,
sieve procedures cannot distinguish between numbers with an even number
of prime factors and an odd number of prime factors, a property
known as the ``parity problem''.  Bombieri (\cite{Bo1}, \cite{Bo2})
clarified things further, showing essentially that knowledge of
$R(\nu)$ for all $\nu<1$ (and no other information about the 
sequence) implies an asymptotic formula for
$\sum_{n\le x} a_n f(n)$ if and only if $f$ gives ``equal weight'' to
numbers with an even number of prime factors and an odd number of prime
factors.
The generalized von Mangoldt functions 
\be\label{Lambdak}
\Lambda_k(n) = \sum_{d|n} \mu(d) \log^k(n/d)
\ee
have this property for $k\ge 2$ (in fact these functions together
with convolutions of the type $\Lambda_{i_1} * \cdots * \Lambda_{i_j}$
($i_1+\cdots+i_j \ge 2$)
form a kind of basis for all such $f$; see \cite{Bo2} for details).
In particular, Bombieri proved that if $R(\nu)$ for all $\nu<1$, then
\be\label{Sk}
S_k(x) := \sum_{n\le x} a_n \Lambda_k(n)
\sim k H A(x) (\log x)^{k-1}.
\ee
A different proof of \eqref{Sk} was given by Friedlander and Iwaniec 
\cite{FI1}.  The required conditions on $A(x)$ and $g$ differ in
\cite{Bo2} and \cite{FI1}, but they are all trivially satisfied if
$A(x)=x$ and $g(d)=1/d$ (here $H=1$).

The special case of \eqref{Sk} corresponding to $k=2$ and $a_n=1$ for all $n$
was earlier proved by Selberg, and it served as a foundation
for the first ``elementary'' proofs of the Prime Number Theorem.

It is natural to inquire what may be deduced from $R(\nu)$
for some \emph{fixed} $\nu<1$.  For twin primes, Bombieri \cite{Bo1}
deduced from $R(\nu)$ for $\nu<1/2$ that
$$
1 - c_k \le \frac{S_k(x)}{kHA(x) (\log x)^{k-1}} \le 1 + c_k,
$$
where $c_2,\cdots$ are constants with $c_k\to 0$ as $k\to\infty$.
We show that knowing $R(\nu)$ for any fixed $\nu<1$ is not
sufficient to deduce \eqref{Sk} for any $k$. 

%
%

\begin{theorem} \label{thm1}
Fix $\nu\in (0,1)$.  There is a sequence $(a_n)$ which satisfies
$R(\nu)$ with $A(x)=x$ and $g(d)=1/d$,
and for which \eqref{Sk} fails for every $k\ge 1$.
Furthermore, we can specify the manner in which
\eqref{Sk} fails, constructing $(a_n)$ so that
\benn
T_k(x) := \frac{S_k(x)}{kx(\log x)^{k-1}}
\eenn
satisfies either (i) $T_k(x) \sim \xi_k$ with $\xi_k<1$ for every $k$; or
(ii) $T_k(x) \sim \xi_k$ with $\xi_k>1$ for every $k$; or (iii) for every
$k$, $\ds \limsup_{x\to\infty} T_k(x)>1$ and $\ds \liminf_{x\to\infty} T_k(x)
<1$.
\end{theorem}

By slightly modifying the construction of the sequence
$(a_n)$, we can create sequences satisfying Theorem 1 for which
$a_n \in \{0,1,2\}$ for every $n$.

Recently there was a major breakthrough on the parity problem by
Friedlander and Iwaniec \cite{FI2}.  
They proved $S_1(x)\sim HA(x)$ under two major assumptions.
First, $R(\nu)$ holds for some $\nu>2/3$.
Second, the bilinear sum condition
\be\label{bilinear}
\sum_{m} \biggl| \sum_{\substack{N<n\le 2N \\ mn\le x }}
\gamma(n,C) \mu(mn) a_{mn} \biggr| \ll A(x) (\log x)^{-1996}, \quad
\gamma(n,C)=\sum_{d|n,d\le C} \mu(d)
\ee
holds uniformly for $\Delta^{-1} x^{\nu/2} < N < \delta^{-1} \sqrt{x}$,
$1\le C \le x^{1-\nu}$, where $\delta, \Delta$ are parameters depending on
 $x$ in such a way that $\delta \to \infty$
and  $\frac{\log \delta}{\log \Delta} \to 0$ 
as $x\to\infty$.  In \cite{FI3}, they applied this successfully to
give an asymptotic formula for the number of primes of the form $a^2+b^4$
which are $\le x$.  The condition \eqref{bilinear} strongly eliminates
the possibility of the sequence having a ``parity bias'', meaning a
tendency for $\mu(n) a_n$ to be of one sign.

The sequences used to prove Theorem \ref{thm1} all exhibit a ``global parity
bias'', meaning that
\be\label{mu}
P(x) = \sum_{n\le x} a_n \mu(n) 
\ee
is large (or large infinitely often).
In light of Selberg's example and the theorem of Friedlander and Iwaniec,
it is natural to inquire whether or not, for each $\nu<1$,
there are sequences $(a_n)$ satisfying $R(\nu)$ and also
\be\label{Pcond}
P(x) \ll_B x(\log x)^{-B} \quad (\forall B>0),
\ee
but failing \eqref{Sk}.
We cannot as yet answer this question entirely, but for all
$\nu<1$, we can construct sequences satisfying $R(\nu)$ and \eqref{Pcond},
but failing \eqref{Sk} for all $k\ge 2$.  These sequences
do satisfy \eqref{Sk} for $k=1$.

\begin{theorem}\label{thm2}
Fix $\nu\in(0,1)$.  There is a sequence $(a_n)$ which satisfies 
$R(\nu)$, \eqref{Pcond} and for which \eqref{Sk} fails for all $k\ge 2$.
\end{theorem}

It is an interesting problem to examine the situation if \eqref{Pcond}
is replaced by a stronger condition (but one weaker than \eqref{bilinear}).
One possibility, suggested by C. Hooley, is to postulate
that the parity bias in arithmetic progressions is small on average,
something like
$$
\sum_{d\le x^{\alpha}} \biggl| \sum_{\substack{n\le x \\ d|n
 }} \mu(n) a_n \biggr| \ll_B x(\log x)^{-B} \quad (\forall B>0).
$$
The sequences we construct for the proof of Theorem \ref{thm1} do
satisfy this condition with arbitrary but fixed $\alpha < 1-\nu$.
The case $\alpha + \nu > 1$ remains open.

%
%
\section{Overall plan}
%
%

The only analytic tool we require is the Prime Number Theorem with 
the de la Vall\'ee Poussin error term.  In fact a much weaker
error term would suffice.

\begin{lemma} \label{PNT}
For some positive constant $c_0$,
\benn
\sum_{n\le x} \Lambda(n) = x+O\( x e^{-c_0\sqrt{\log x}}\).
\eenn
\end{lemma}

Assume without loss of generality that $\nu>1/2$.
Let $M$ be an integer, and $\del$ and
$\varpi$ be real numbers satisfying
\be
0 < \del \le \frac{1}{3M^2}, \quad M\del+\frac{1}{M} < \varpi < 1-\nu.
\label{eq:Mdelpi}
\ee
Take $x_0$ sufficiently large and 
$c_1 \in (0,c_0)$ (both depending on $\nu$, $M$ and $\del$).  For $j\ge 1$ put
\be
x_{j+1}=x_j \( 1 + e^{-c_1\sql{x_j}}\),
 \quad I_j = \mathbb Z \cap (x_j,x_{j+1}],\quad K_j = |I_j|.
\label{eq:xIH}
\ee
In what follows, all constants implied by the $O-$  symbol may
depend on $\nu$ and $M$.  Dependence on other variables will
be indicated by subscripts to the $O-$ symbol.  The numbers
$a_n$ for $n\in I_j$ will satisfy three basic properties.  First,
\be
0 \le a_n \le 2.
\label{eq:an_bounded}
\ee
Second,
\be
\sum_{\substack{n\in I_j \\ d|n}} a_n = \frac{K_j}{d} + O\( \frac{K_j}{d}
e^{-c_1\sql{x_j}} \) \qquad (1 \le d\le x_j^{1-\varpi}).
\label{eq:an_d|n}
\ee
Third, for some positive constants $\theta_k$ ($k\ge 1$ for Theorem 
\ref{thm1}, $k\ge 2$ for Theorem \ref{thm2}) which depend
on $\nu$, $M$ and $\del$, and some numbers $\sigma_j \in \{-1, 1\}$ (which we 
are free to choose), we have
\be
\sum_{n\in I_j} a_n \Lambda_k(n) = (k+\sigma_j\theta_k) K_j (\log x_j)^{k-1}
 \(1 + O_k \( e^{-c_1\sql{x_j}} \) \).
\label{eq:sum_in_Ij}
\ee
For Theorem 2, we also require that
\be
\sum_{n\in I_j} a_n \mu(n) = O\( K_j e^{-c_1\sql{x_j}}\).
\label{eq:anmun_Ij}
\ee
Deducing Theorems \ref{thm1} and \ref{thm2} from \eqref{eq:an_bounded}--
\eqref{eq:anmun_Ij} is straightforward.  For
$d\le x^{\nu}$, \eqref{eq:an_bounded} implies
\benn
A_d(x) = O\( \frac{x^{\nu+\varpi}}{d}\) + \sum_{x^{\nu+\varpi}\le x_j
\le x}\;\; \sum_{\substack{n\in I_j \\ d|n}} a_n
\eenn
By \eqref{eq:Mdelpi}, if $x_j\ge x^{\nu+\varpi}$, then 
$d < x^{(1-\varpi)(\nu+\varpi)} \le x_j^{1-\varpi}$.
Thus, by \eqref{eq:an_d|n},
\benn
\begin{split}
A_d(x) &=  O\( \frac{x^{\nu+\varpi}}{d}\)  + \sum_j \frac{K_j}{d} \( 1 +
  O\(e^{-c_1\sqrt{(\nu+\varpi)\log x}}\)\) \\
&= \frac{x}{d} + O\( \frac{x}{d} e^{-\frac12 c_1\sql{x}} \) = 
 g(d) A(x) + O\( \frac{A(x)}{d} e^{-\frac12 c_1\sql{x}} \).
\end{split}
\eenn
Summing on $d$ gives $R(\nu)$.  Similarly,  \eqref{eq:anmun_Ij} implies
\eqref{Pcond}.  
From \eqref{Lambdak}, we have $\log^k n =(1*\Lambda_k)(n) \ge \Lambda_k(n)$.
Thus, using \eqref{eq:an_bounded} and \eqref{eq:sum_in_Ij}, we obtain
\benn
\begin{split}
\sum_{n\le x} a_n \Lambda_k(n) &= O(x(\log x)^{k-2})
  +\sum_{\frac{x}{\log^2 x}\le x_j\le x}\;\; \sum_{n\in I_j} a_n \Lambda_k(n)\\
&=  O(x(\log x)^{k-2}) + (\log x)^{k-1} \(1 +
  O\(\frac{\log\log x}{\log x}\)\) \sum_{x_j\le x } (k+\sigma_j\theta_k) K_j\\
&= (\log x)^{k-1} (xk + \theta_k \sum_{x_j\le x} \sigma_j K_j) +
  O(x(\log x)^{k-3/2}).
\end{split}
\eenn
The three types of behavior for $T_k(x)$ in Theorem 1
are obtained by taking (respectively)
(i) $\sigma_j=-1$ for all $j$; (ii) $\sigma_j=1$ for all $j$; or (iii)
$\sigma_j=-1$ if $2^{2^r} < x_j \le 2^{2^{r+1}}$ for an even $r$ and
$\sigma_j=1$ if $2^{2^r} < x_j \le 2^{2^{r+1}}$ for an odd $r$.

It remains, therefore, to construct numbers
$a_n$ on each interval
$I_j$ satisfying \eqref{eq:an_bounded}--\eqref{eq:anmun_Ij} as appropriate for
Theorems \ref{thm1} and \ref{thm2}.
The basic idea is to start with $a_n=1$ for all $n$, then
shift around some of the mass from the numbers $a_n$ with $n$ composed of
``large'' prime factors.  This must be done
very delicately in order to preserve \eqref{eq:an_d|n}, and this is
the most complex part of the argument.
 We will
work with smooth functions defined on numbers with a given number of
prime factors.  Let
\begin{gather*}
T_r = \{ (u_1,\ldots,u_r): 0\le u_1 \le \cdots \le u_r, u_1+\cdots+u_r=1\},\\
U_r = \{ (u_1,\ldots,u_r): u_i \ge 0\; (1\le i\le r), u_1+\cdots+u_r=1\}.
\end{gather*}

For positive numbers $\e, B$, let $\FF_r(\e,B)$ be the set of functions
$f(u_1,\ldots,u_r)$ on $U_r$ that are (i) symmetric in all variables,
 (ii) zero whenever $\min u_i \le \e$ and (iii) $f$ and all first order partial
derivatives are at most $B$ in absolute value on $U_r$.
If $n=p_1 \cdots p_r$, the numbers $p_i$ being primes with
no assumptions on their relative sizes, then
\benn
f \( \frac{\log p_1}{\log n}, \ldots, \frac{\log p_r}{\log n} \)
\eenn
is well-defined.  With these assumptions, we may estimate in a standard way
sums over $f$ in terms of integrals.   

%
%

\begin{lemma}\label{lem:sum f}  Let $f\in \curly{F}_r(\e,B)$, $0\le y\le x$
and $x$ large in terms of $\e,r,B$.  Then
\benn
\begin{split}
\sum_{\substack{p_1,\cdots,p_r \\ x\le n=p_1\cdots p_r \le x+y}} & f \(
   \frac{\log p_1}{\log n}, \ldots, \frac{\log p_r}{\log n} \)  \\
&= \frac{y}{\log x} \int_{U_r} \frac{f(u_1,\ldots,u_r)}
  {u_1\cdots u_r} + O_{\e,r,B} \( \tfrac{y^2}{x\log x} + x 
  e^{-c_0\sqrt{\e^{r-1}\log x}} \).
\end{split}
\eenn
When $r=1$, the integral is $f(1)$.
\end{lemma} 

\begin{proof}  Let $F$ denote the sum in the lemma.
In this proof, constants  implied by the $O-$ symbol may depend
on $\e,r,B$.  When $r=1$, by Lemma \ref{PNT},
$$
F = f(1) (\pi(x+y)-\pi(x)) = \frac{y}{\log x} f(1) + O(xe^{-c_0\sql{x}}).
$$
We now proceed by induction on $r$.  Suppose $r\ge 2$ and fix 
$p_1 \in [x^{\e}, 2x^{1-\e}]$.  
Writing $n'=p_2\cdots p_r$ and $v_j=\frac{\log p_j}
{\log n'}$ for $2\le j\le r$, we have
$$
f\(\tfrac{\log p_1}{\log n}, \cdots, \tfrac{\log p_r}{\log n}\) 
=\(1+ O\pfrac{y}{x\log x}\) g(v_2,\cdots,v_r),
$$
where
$$
g(v_2,\cdots,v_r)=f\(\tfrac{\log p_1}{\log x}, \tfrac{\log(x/p_1)}{\log x} v_2,
\ldots \tfrac{\log(x/p_1)}{\log x} v_r \).
$$
We have $g\in \FF_{r-1}(\e,B)$, so by the induction hypothesis 
$$
F = \sum_{p_1} \frac{y}{p_1 \log(x/p_1)} \int\limits_{\vv \in U_{r-1}} 
  \frac{g(v_2,\ldots,v_r)}{v_2 \cdots v_r} + O\( \frac{y^2}{p_1 x\log x} + 
  \frac{x}{p_1}e^{-c_0 \sqrt{\e^{r-2}\log(x/p_1)}} \).
$$
Since $\log(x/p_1) \ge \e \log x -1$, $\sum_{p_1} 1/p_1 \ll 1$ and thus
the error terms above total 
$$
O\( \frac{y^2}{x\log x} + xe^{-c_0 \sqrt{\e^{r-1}\log x}} \).
$$
By Lemma \ref{PNT} and partial summation, for a fixed $v_2,\ldots,v_r$, we have
\begin{align*}
\sum_{p_1} &\frac{f(\tfrac{\log p_1}{\log x}, \tfrac{\log(x/p_1)}{\log x} v_2,
  \ldots, \tfrac{\log(x/p_1)}{\log x} v_r)}{p_1 \log(x/p_1)} \\
&\qquad = \int_{x^{\e}}^{2x^{1-\e}} \frac{f(\tfrac{\log t}{\log x}, 
  \tfrac{\log(x/t)}{\log x} v_2,\ldots, \tfrac{\log(x/t)}{\log x} v_r)}
  {t \log t\log(x/t)}\, dt + O(e^{-c_0\sqrt{\e \log x}}) \\
&\qquad = \frac{1}{\log x} \int_{\e}^{1-\e/2} \frac{f(u,(1-u)v_2,\ldots,
  (1-u)v_r)}{u(1-u)}\, du +  O(e^{-c_0\sqrt{\e \log x}}).
\end{align*}
Therefore
$$
F = \frac{y}{\log x} \int\limits_{\substack{v_2+\cdots+v_r=1 \\ 0<u<1}} 
  \frac{f(u,(1-u)v_2,\ldots,(1-u)v_r)}{u(1-u)v_2\cdots v_r} 
  + O\( \frac{y^2}{x\log x} + xe^{-c_0\sqrt{\e^{r-1}\log x}} \).
$$
Making the change of variables $u_1=u$, $u_j=(1-u)v_j$
($2\le j\le r$) gives the lemma.
\end{proof}

%
%
%
%
\section{The construction on $I_j$}
%
%
%
%

To facilitate working with sets of numbers with prime factors in specific
ranges, we adopt some special
 notation.  A \textit{partition} is a non-decreasing
sequence of positive integers $\ba=(\a_1,\ldots,\a_r)$ (also thought of as
a ``multi-set'').  Let $|\ba| = r$ and $\Sigma(\ba)=\a_1+\cdots+\a_r$.
Let $\perm(\ba)$ be the number of permutations of the numbers in $\ba$,
e.g. $\perm(1,1,2,3)=12$.  Let 
\benn
\one_n = \overbrace{(1,\ldots,1)}^n
\eenn
and let $P(m)$ be the set of all partitions of $m$ (all $\ba$ with
$\Sigma(\ba)=m$).  Let $\bE$ denote the empty partition ($|\bE|=0$
and $\Sigma(\bE)=0$).
The notation $\ba \subseteq \bb$ means that each number
in $\bb$ occurs at least as many times as the number occurs in $\ba$,
and $\ba+\bb$ is the partition consisting of all the parts of $\ba$ and
of $\bb$, so in particular $|\ba+\bb|=|\ba|+|\bb|$ and
$\Sigma(\ba+\bb)=\Sigma(\ba)+\Sigma(\bb)$.  Also, if $\ba \subseteq \bb$,
$\bb-\ba$ is defined by $\ba + (\bb-\ba)=\bb$.

For brevity, write $x=x_j$, $K=K_j$, $I=I_j$.  For $1\le i\le M$, let
$\curly{P}_i$ be the set of primes in the interval
$[x^{i(1/M-\del)}, x^{i(1/M+\del)}]$.  For each
partition $\ba=(\a_1,\cdots,\a_r)$, let
\benn
\curly{D}_\ba = \{ p_1\cdots p_r :  p_i \in \curly{P}_{\a_i} \;
(1\le i\le r) \}, \quad \curly{C}_\ba = \curly{D}_\ba \cap I.
\eenn
In particular, $\curly{D}_\bE = \{ 1 \}$ and $\curly{C}_{(M)}$ is the set 
of primes in $I$.  Also, by
\eqref{eq:Mdelpi}, $\curly{C}_\ba$ is empty unless $\ba \in P(M)$.

Let $c_1 = c_0 (2M)^{-r}$.
We put $a_n=1+b_n$, where $|b_n| \le 1$, and $b_n=0$ unless 
$n$ lies in some $\curly{C}_\ba$ with $\ba\in P(M)$.
Thus, if $1<d\le x^{1-\varpi}$, then
\benn
\sum_{\substack{n\in I \\ d|n}} a_n = \frac{K}{d}+O(1)
\eenn
unless $d\in\curly{D}_\bb$ for some $\bb=(\b_1,\ldots,\b_s)$.
In this case $x^{(1/M-\del)(\b_1+\cdots+\b_s)} \le d \le x^{1-\varpi}$,
which by \eqref{eq:Mdelpi} implies $\b_1 + \cdots + \b_s \le M-2$. 
Let $Q=\bE \cup P(1) \cup \cdots \cup P(M-2)$.  To obtain 
\eqref{eq:an_d|n}, it suffices to prove that for each $\bb\in Q$ and each $d\in
\curly{D}_\bb$,
\be
\sum_{\substack{\ba \in P(M) \\ \bb \subseteq \ba}} \;\; \sum_{\substack{
n\in \curly{C}_\ba \\ d|n}} b_n = O \( \frac{K}{d} e^{-c_1\sql{x}} \).
\label{eq:mainsumbn}
\ee
This system of inequalities has the trivial solution $b_n=0$ for all $n$,
but we need a solution with $|b_n| \gg 1$ on average in order to obtain
\eqref{eq:sum_in_Ij}.

For $1\le i\le M$, let $J_i=[i(1/M-\del),i(1/M+\del)]$.  For each
$\ba=(\a_1,\ldots,\a_r)\in P(M)$, suppose $f_\ba\in 
\curly{F}_r(\tfrac{1}{2M},B)$ is
supported on $T_r \cap (J_{\a_1} \times \cdots \times J_{\a_r})$ and
the symmetric regions in $U_r$.  For $n\in \curly{C}_\ba$,
$n=p_1\cdots p_r$, let
\be
b_n = f_\ba \( \frac{\log p_1}{\log n}, \cdots , \frac{\log p_r}{\log n}\).
\label{eq:bndef}
\ee
Suppose that $\bb=(\b_1,\ldots,\b_s)\in Q$ with $\bb \subseteq \ba$.
Then $r\ge s+1$. Let
$d=p_1 \cdots p_s \in \curly{D}_\bb$ with $p_i\in \curly{P}_{\b_i}$
$(1\le i\le s)$ and put $v_i=\frac{\log p_i}{\log x}$ for $1\le i\le s$.
We have
$$
f_\ba \( \tfrac{\log p_1}{\log n}, \cdots , \tfrac{\log p_r}{\log n}\) =
g\( \tfrac{\log p_{s+1}}{\log (n/d)}, \ldots, \tfrac{\log p_{r}}{\log (n/d)}\)
+ O(e^{-c_1\sqrt{\log x}}),
$$
where
$$
g(w_1,\ldots,w_{r-s}) = f(v_1,\cdots,v_s,\tfrac{\log(x/d)}{\log x} w_1,
\ldots,\tfrac{\log(x/d)}{\log x} w_{r-s}).
$$
Since $g\in \FF_{r-s}(\tfrac{1}{2M},B)$, Lemma \ref{lem:sum f} implies that
\begin{align*}
\sum_{\substack{n\in \curly{C}_{\ba} \\ d|n}} b_n &=
  \frac{K}{(r-s)! d \log(x/d)} \int\limits_{\ww \in U_{r-s}} 
  \frac{g(w_1,\ldots,w_{r-s})}{w_1 \cdots w_{r-s}} + 
  O\( \frac{x}{d} e^{-2c_1\sqrt{\log x}} \) \\
&= \frac{K}{(r-s)! d\log x} \int_{\uu \in V_{r-s}(1-v_1-\cdots-v_s)} 
  \frac{f_{\ba} (v_1,\ldots,v_s,u_1,\cdots,u_{r-s})}{u_1 \cdots u_{r-s}}
  + O\( \frac{K}{d} e^{-c_1\sql{x}}\),
\end{align*}
where $V_t(A)=\{ (u_1,\ldots,u_t): u_i\ge 0 \forall i, \sum u_i = A \}$.

Therefore, to prove \eqref{eq:mainsumbn}, it suffices to find functions
$f_\ba$ so that for all $\bb=(\b_1,\ldots,\b_s) \in Q$ and
$(v_1,\ldots,v_s)\in J_{\b_1} \times \cdots \times J_{\b_s}$, we have
\be
\sum_{\bmu \in P(M-\Sigma(\bb))}
  \frac{1}{|\bmu|!} \int\limits_{\uu \in V_{|\bmu|}(1-v_1-\cdots-v_s)}
  \frac{f_{\bb+\bmu}(v_1,\ldots,v_s,u_1,\ldots,u_{r-s})}{u_1\cdots u_{r-s}}
   = 0.
\label{eq:mainint}
\ee
When $M\ge 6$, $|Q| > |P(M)|$ (i.e. there
are more equations than functions), but there
is enough structure in the system \eqref{eq:mainint} to find a nontrivial
solution.   In fact, once $f_{\one_M}$ is chosen, the other functions
$f_\ba$ are uniquely determined by \eqref{eq:mainint}, but we do not need to
prove this.
Suppose $\ba=(\a_1,\cdots,\a_r) \in P(M)$, $\brho$ is a permutation of
$\ba$ and $\vv \in J_{\rho_1} \times \cdots \times J_{\rho_r}$.  For
some constant $e_\ba$, define
\be
f_\ba(\vv) = e_\ba v_{1} \cdots v_r \int\limits_{\eqref{wsums}}
\frac{f_{\one_M}(\mathbf{w})}{\prod w_{ij}},
\label{eq:fdef}
\ee
where the integration is over the set of  $\ww = \{ w_{ij}: 
1\le j\le r,1\le i\le \rho_j \} \in J_1^M$ with
\be\label{wsums}
\sum_{i=1}^{\rho_j} w_{ij} = v_j \quad (1\le j\le r).
\ee
For example, if $\ba=(1,1,2,3)$, $v_1\in J_1$, $v_2\in J_2$, $v_3\in J_3$
and $v_4\in J_1$, we have
\benn
f_\ba(\vv) = e_\ba v_1 v_2 v_3 v_4 \int\limits_{\substack{w_{11}=v_1 \\
w_{12}+w_{22}=v_2 \\ w_{13}+w_{23}+w_{33}=v_3 \\ w_{14} = v_4}}
\frac{f_{\one_7}(w_{11}, w_{12}, w_{22}, w_{13}, w_{23}, w_{33}, w_{14})}
{w_{11} w_{12} w_{22} w_{13} w_{23} w_{33} w_{14} }.
\eenn
For consistency, set $e_{\one_M}=1$.

We next show that
substituting \eqref{eq:fdef} into \eqref{eq:mainint} reduces the
problem to solving a system of equations in the numbers $e_\ba$.  
Fix $\bb=(\b_{1},\ldots,\b_s) \in Q$ and $\bmu=(\mu_{1},\ldots,
\mu_{r-s}) \in P(M-\Sigma(\bb))$.
 Suppose that $v_j \in J_{\b_j}$ ($1\le j\le s$) and let $\brho$ be a
permutation of $\bmu$.
Take $\uu$ so that $\sum v_i + \sum u_i = 1$
and $u_i \in J_{\rho_i}$ ($1\le i\le r-s$).
Because $f_{\bb+\bmu}$ is symmetric in all
variables, for each $\brho$ the contribution to the integral
in \eqref{eq:mainint} is identical.  In other words, the integral
in  \eqref{eq:mainint} equals $\perm(\bmu)$ times the integral over
those $\uu\in J_{\mu_1} \times \cdots \times J_{\mu_{r-s}}$.
For such $\uu$, \eqref{eq:fdef} implies
\benn
f_{\bb+\bmu}(\vv,\uu)=e_{\bb+\bmu} v_{1}\cdots v_s u_{1} \cdots 
  u_{r-s} \int\limits_{\eqref{eq:intcond1}} \frac{f_{\one_M}(\ww,\zz)}
  {\prod_{i,j} w_{ij} \prod_{i,j} z_{ij}},
\eenn
where the integral is over the variables $w_{ij}$, $z_{ij} \in J_1$ satisfying
\be
\sum_{i=1}^{\b_j} w_{ij}=v_j \;\; (1\le j\le s); \quad
\sum_{i=1}^{\mu_j} z_{ij}=u_j \;\; (1\le j\le r-s).
\label{eq:intcond1}
\ee
Thus, with $\bb$ and $\vv$ fixed,
$$
\int\limits_{\substack{\uu \in V_{|\bmu|}(1-v_1-\cdots-v_s) \\ u_j \in 
  J_{\mu_j} (1\le j\le r-s)}}  \!\!\!\! \frac{f_{\bb+\bmu} (\vv,\uu)}
  {u_1\cdots u_{r-s}}
= e_{\bb+\bmu} v_{1} \cdots v_s \int\limits_{\substack{
  \uu \in V_{r-s}(1-v_1-\cdots-v_s) \\ \eqref{eq:intcond1}}}
  \frac{f_{\one_M} (\ww,\zz)} {\prod w_{ij} \prod z_{ij}}.
$$
Since $u_1,\ldots,u_{r-s}$ are dependent variables in the integral 
on the right side, the left side is actually independent of $\bmu$.
Thus, \eqref{eq:mainint} follows from the system
\be
\sum_{\bmu\in P(M-\Sigma(\bb))} \frac{\perm(\bmu)}{|\bmu|!} e_{\bb+\bmu} 
= 0 \quad (\bb\in Q), \quad e_{\one_M}=1.
\label{eq:ealpha_equn}
\ee
As noted before, \eqref{eq:ealpha_equn} has more equations than variables
when $M>6$, but there is a simple solution 
(again we do not need to prove uniqueness, but it is straightforward),
namely 
\be
e_\ba = \frac{(-1)^{\Sigma(\ba)+|\ba|}}{\a_1 \cdots \a_r}, \qquad
\ba=(\a_1,\ldots,\a_r).
\label{eq:ealpha def}
\ee
With \eqref{eq:ealpha def}, $e_{\bb+\bmu}=e_\bb e_\bmu$  for all
$\bb,\bmu$, so \eqref{eq:ealpha_equn} is equivalent to 
$\gamma_m=0$ $(2\le m\le M)$, where
\benn
\gamma_m := \sum_{\substack{\bmu\in P(m) \\ \bmu=(\mu_1,\ldots,\mu_r)}}
\frac{(-1)^{|\bmu|}}{|\bmu|!} \frac{\perm(\bmu)}{\mu_1\cdots \mu_r} =
\sum_{r=1}^m \frac{(-1)^r}{r!} \sum_{\substack{d_1+\cdots+d_r=m\\
d_i \ge 1 (1\le i\le r)}} \frac{1}{d_1\cdots d_r}.
\eenn
This follows by considering the generating function
\benn
G(z) = \sum_{m=1}^\infty \gamma_m z^m.
\eenn
Since $|\gamma_m| \le \sum_{r=1}^m \binom{m-1}{r-1} = 2^{m-1}$, $G(z)$
has radius of convergence $\ge 1/2$.  Thus, for $|z| \le 1/3$,
\benn
\begin{split}
G(z) &= \sum_{r=1}^\infty \frac{(-1)^r}{r!} \sum_{d_1,\ldots,d_r \ge 1} 
\frac{z^{d_1+\cdots + d_r}}{d_1\cdots d_r} \\
&=  \sum_{r=1}^\infty \frac{(-1)^r}{r!} \( \sum_{d=1}^\infty \frac{z^d}{d} \)
^r \\
&= \sum_{r=1}^\infty \frac{(-1)^r}{r!} \( -\log(1-z) \)^r = e^{\log(1-z)}-1
=-z,
\end{split}
\eenn
which proves \eqref{eq:ealpha_equn}.  As noted earlier, \eqref{eq:ealpha_equn}
implies \eqref{eq:mainint}, which implies \eqref{eq:mainsumbn}, which
implies \eqref{eq:an_d|n}.

Modulo the choice of function $f_{\one_M}$, we have constructed our
numbers $b_n$.  The following theorem sums up the properties we are
interested in.

%
%

\begin{theorem}\label{thm3}  Fix $M,\varpi,\del$ so that
\eqref{eq:Mdelpi} is satisfied and also $\del \le (2M)^{-M}$.
Let $B$ be large depending on $M,\del$.
Let $f_{\one_M} \in \curly{F}_M(\frac{1}{2M},B)$ 
with $|f_{\one_M}(\uu)| \le 1$ for all
$\uu \in U_M$.  For every $\ba \in P(M)$, define $e_\ba$ by
\eqref{eq:ealpha def}, define $f_\ba$ by \eqref{eq:fdef},
$b_n$ by \eqref{eq:bndef}, and put $a_n=1+b_n$.  Then, for each 
interval $I=I_j$, 
\eqref{eq:an_bounded} and \eqref{eq:an_d|n} are satisfied, plus we have
\be
\sum_{n\in I} a_n \Lambda_k(n) = K(\log x)^{k-1} \bigl[ k + (-1)^{M+1}Z_k +
O(e^{-c_1\sql{x}}) \bigr] \qquad (k\ge 1)
\label{eq:sum anLk}
\ee
and
\be
\sum_{n\in I} a_n \mu(n) = \frac{K}{\log x}\bigl[(-1)^M Z_0 +
O(e^{-c_1\sql{x}}) \bigr],
\label{eq:sum anmun}
\ee
where
\benn
Z_k := \int_{\uu \in U_M}\frac{u_1^k f_{\one_M}(\uu)}{u_1\cdots u_M} 
 = \frac{1}{M} \int_{\uu \in U_M} \frac{u_1^k + \cdots + u_M^k}
{u_1\cdots u_M} f_{\one_M}(\uu). 
\eenn
Formula \eqref{eq:sum anmun} also holds with $\mu(n)$
replaced by $\lambda(n)$.
\end{theorem}

\begin{proof}  We have already seen that \eqref{eq:an_d|n} is satisfied.
Let $\ba \in P(M)$, $\ba \ne \one_M$.
By \eqref{eq:ealpha def},  $|e_\ba| \le 1$, so by \eqref{eq:fdef},
\benn
|f_\ba| \le (1/M-\del)^{-M} \del^{M-|\ba|}
\le \del (2M)^M \le 1.
\eenn
Next, by \eqref{Lambdak} and Lemma \ref{lem:sum f}, 
for each $\ba=(\a_1,\cdots,\a_r)\in P(M)$, we have
\benn
\begin{split}
\sum_{n\in \curly{C}_\ba} b_n \Lambda_k(n) &= \frac{1}{r!} \sum_{\substack{
  p_1,\cdots,p_r \\ n=p_1 \cdots p_r \in \curly{C}_\ba}} f_\ba\( 
  \tfrac{\log p_1}{\log n}, \cdots, \tfrac{\log p_r}{\log n}\) 
  \sum_{\eps_1,\cdots,\eps_r \in \{0,1\}} (-1)^{r-\e_1-\cdots-\e_r} 
  \log^k (p_1^{\e_1} \cdots p_r^{\e_r}) \\
&= \frac{(-1)^r\perm(\ba)}{r!} K(\log x)^{k-1} \biggl[ O(e^{-c_1\sql{x}}) \\
  &\qquad + \!\!\int\limits_{\substack{\vv \in U_r \\ v_i\in J_{\a_i}\; 
  (1\le i\le r)}} \frac{f_{\ba}(\vv)}{v_1\cdots v_r} \sum_{\eps_1,\cdots,
  \eps_r \in \{0,1\}} (-1)^{\eps_1+\cdots+\eps_r} \biggl( \sum_{j=1}^r
  \eps_j v_j \biggr)^k \biggr].
\end{split}
\eenn
By \eqref{eq:fdef} and the fact that $f_{\one_M}$ is symmetric in all
variables, we obtain
\benn
\begin{split}
\sum_{n\in \curly{C}_\ba} b_n \Lambda_k(n) &=
  \frac{\perm(\ba) e_\ba (-1)^{|\ba|}}{|\ba|!} K
  (\log x)^{k-1} \biggl[ O(e^{-c_1\sql{x}}) \\
&\qquad + \int\limits_{\ww \in U_M}
  \!\! \frac{f_{\one_M}(\ww)}{\prod w_{i,j}} \sum_{\boldsymbol{\eps}}
  (-1)^{\eps_1+\cdots+\eps_r} \biggl( \sum_{j=1}^{|\ba|} \eps_j 
  (w_{j,1}+\cdots+w_{j,\a_j}) \biggr)^k \biggr] \\
&= \frac{\perm(\ba) e_\ba (-1)^{|\ba|}}{|\ba|!} K (\log x)^{k-1} \biggl[
  O(e^{-c_1\sql{x}}) \\
&\qquad + \int\limits_{\uu\in U_M} \frac{f_{\one_M}(\uu)}{u_1\cdots u_M} 
  \sum_{N=1}^M \biggl( \!\!
 \sum_{\substack{\eps_1,\ldots,\e_r \in \{0,1\} \\ \eps_1\a_1+\cdots+\eps_r\a_r
  =N}}\!\! (-1)^{\eps_1+\cdots+\eps_r} \biggr) \( u_1 + \cdots + u_N \)^k
  \biggr].
\end{split}
\eenn
Summing on $\ba\in P(M)$ and using \eqref{eq:ealpha def} gives
\begin{multline}
\sum_{n\in I} b_n \Lambda_k(n) = K (\log x)^{k-1} \biggl[ O(e^{-c_1\sql{x}})\\
  + (-1)^M \sum_{N=1}^M W(M,N) \int\limits_{\uu\in U_M}
  \frac{f_{\one_M}(\uu)}{u_1\cdots u_M} (u_1+\cdots+u_N)^k 
  \biggr],
\label{eq:sum bnLk}
\end{multline}
where
\be
W(M,N) = \sum_{r=1}^M \frac{1}{r!} \sum_{\substack{d_1+\cdots+d_r=M \\
  d_i\ge 1\; \forall i}} \frac{1}{d_1\cdots d_r} 
  \sum_{\substack{\eps_1, \ldots, \eps_r\in \{0,1\} \\
  \eps_1 d_1 + \cdots + \eps_r d_r=N}} (-1)^{\eps_1+\cdots +\eps_r}.
\label{eq:W def}
\ee 
By examining the generating function, we next prove that
\be
W(M,N) = 0 \quad (2\le N\le M), \qquad W(M,1)=-1, W(M,0)=1 \quad (M\ge 1).
\label{eq:WMN eq}
\ee
For $\max(|x|,|y|) < 1/3$, we have
\benn
\begin{split}
\sum_{M=1}^\infty \sum_{N=0}^M W(M,N) &x^M y^N = \sum_{r=1}^\infty 
  \frac{1}{r!} \sum_{\eps_1,\ldots,\eps_r\in\{0,1\}}
   (-1)^{\eps_1+\cdots+\eps_r} \sum_{d_1,\ldots,d_r \ge 1}
  \pfrac{x^{d_1} y^{\eps_1 d_1}}{d_1}
  \cdots \pfrac{x^{d_r} y^{\eps_r d_r}}{d_r} \\
&= \sum_{r=1}^\infty \frac{1}{r!} \sum_{\eps_1,\ldots,\eps_r\in \{0,1\}} (-1)^{
  \eps_1+\cdots+\eps_r} \( -\log(1-xy^{\eps_1}) \) \cdots 
 \( -\log(1-xy^{\eps_r}) \) \\
&=  \sum_{r=1}^\infty \frac{1}{r!} \( \log(1-xy) - \log(1-x) \)^r \\
&= -1 + \frac{1-xy}{1-x} = (1-y)(x+x^2+x^3 + \cdots).
\end{split}
\eenn
This proves \eqref{eq:WMN eq}, and, together with \eqref{eq:sum bnLk},
completes the proof of \eqref{eq:sum anLk}.   For the sum of $\mu(n) b_n$,
we obtain a similar expression corresponding to the ``$N=0$'' term.
Thus
$$
\sum_{n\in I} b_n \mu(n) = \frac{(-1)^M K}{\log x} \( W(M,0) I_0
+ O(e^{c_1\sql{x}}) \).
$$
The asymptotic \eqref{eq:sum anmun} now follows from \eqref{eq:WMN eq}.
Lastly, $b_n=0$ if $n$ has a prime factor $< x^{1/M-\del}$.  Hence,
when $b_n \mu(n) \ne b_n \lambda(n)$, $n$ is divisible by the square of
a prime $ \ge x^{1/M-\del}$.  The number of such $n\le x$ is $\ll
x^{1-1/M+\del}$ and this proves the final claim.
\end{proof}

%
%

\begin{proof}[Proof of Theorems 1, 2.]
Define 
$$
\ell(v_1,\ldots,v_M;\xi)=\max(0,\xi^{-4}(\xi^2-v_1^2-\cdots-v_M^2)^2),
$$
which is nonzero only when $|v_i| \le \xi$ for each $i$.
To prove Theorem \ref{thm1}, take in Theorem \ref{thm3}
$$
f_{\one_M}(\uu) =(-1)^{M+1} \sigma_j \ell(u_1-1/M,\ldots,u_M-1/M;\del).
$$
For $\uu\in U_M$,  $u_1+\cdots+u_M=1$ and thus $Z_1 = Z_0/M$.
To prove Theorem \ref{thm2}, we must exhibit
a function $f_{\one_M}$ so that $Z_0=0$ and $Z_k\ne 0$ for $k\ge 2$.
Let $M$ be even and put $\ww=(\frac{1}{M},\cdots,\frac{1}{M})$.  Let $V$
be the set of vectors in $\RR^M$ with exactly $M/2$ components equal to
$\del/2$ and $M/2$ components equal to $-\del/2$.
We will take
$$
f_{\one_M}(\uu) = u_1 \cdots u_M \biggl[ \ell(\uu-\ww;\del^3) -
\binom{M}{M/2}^{-1} \sum_{\vv\in V} \ell(\uu-\ww-\vv;\del^3) \biggr].
$$
Letting 
$$
J = \int_{v_1+\cdots+v_M=0} \ell(\vv;\del^3),
$$
it follows that
\begin{align*}
Z_k &\le M(\tfrac{1}{M}+\del^3)^k J - \bigl[ \tfrac{M}{2} (\tfrac{1}{M}-
\tfrac{\del}{2}-\del^3)^k + \tfrac{M}{2}  (\tfrac{1}{M}+
\tfrac{\del}{2}-\del^3)^k \bigr] J \\
&= \frac{J}{2 M^{k-1}} \bigl[ 2(1+\lam)^k - (1-\eps-\lam)^k-(1+\eps-\lam)^k 
\bigr],
\end{align*}
where $\eps=\frac{\del}{2M}$ and $\lambda=\frac{\del^3}{M}$.
Since $x^k$ has convex derivative for $x>0$, we have
\begin{align*}
(1+\eps-\lam)^k-2(1+\lam)^k+(1-\eps-\lam)^k &\ge
  (\eps-2\lam) k (1+\eps/2)^{k-1} - (\eps+2\lam) k (1+\lam)^{k-1} \\
&\ge k(1+\lam)^{k-2} [ (\eps-2\lam)(1+\eps/2)-(\eps+2\lam)(1+\lam) ] \\
&= k(1+\lam)^{k-1} \( \frac{\del^2}{8M^2} + O\pfrac{\del^3}{M} \).
\end{align*}
This proves $Z_k<0$ for $k\ge 2$ if $\del$ is small enough,
and completes the proof.
\end{proof}

\begin{Acknowledgement}
The author thanks John Friedlander for helpful conversations and Denka
Kutzarova for help constructing the function $f_{\one_M}$ for Theorem 2.
\end{Acknowledgement}

%
%
%

\end{document}